\documentclass[draft,12pt]{article}

\usepackage{amsfonts,amsmath,amsthm,amscd,amssymb,latexsym,cite,verbatim,texdraw,floatflt,caption2,pb-diagram}

\usepackage[mathscr]{eucal}
\newtheorem{theorem}{Theorem}[section]
\newtheorem{lemma}{Lemma}[section]
\newtheorem{remark}{Remark}[section]
\newtheorem{corollary}{Corollary}[section]
\theoremstyle{definition}

\newcommand{\keywords}{\textbf{Keywords: } }

\newcommand{\subjclass}{\textbf{Mathematics Subject Classification (2010):} }
\renewcommand{\abstract}{\textbf{Abstract.} }

\numberwithin{equation}{section}

\begin{document}

%%%%%%%%%%%%%%%%%%%%%%%%%%%%%%%%%%%%%%%%%%%%%%%%%%%%%%%%%%%%%%%%%%%%%%%%%%%%%

\title{Approximation  on hexagonal domains by Taylor-Abel-Poisson means}

\author{J\"{u}rgen Prestin, Viktor Savchuk, Andrii Shidlich}

\maketitle

\begin{abstract}
Approximative properties of the Taylor-Abel-Poisson linear summation me\-thod of Fourier series are   considered  for functions of several variables, periodic with respect to the hexagonal domain, in the integral metric. In particular,  direct and inverse theorems are proved in terms of approximations of functions by the Taylor-Abel-Poisson means and  $K$-functionals generated by radial derivatives. Bernstein type inequalities for $L_1$-norm of high-order radial derivatives  of the Poisson kernel are also obtained.

\end{abstract}



%%%%%%%%%%%%%%%%%%%%%%%%%%%%%%%%%%%%%%%%%%%%%%%%%%%%%%%%%%%%%%%%%%%%%%%%%%%%%%%%%%%%%%%%%%%%%%%%%%%%%%%%%%%%%%%%%%%%%%%%%%%%%%%%

\section{Introduction} In this paper, the approximative properties of the Taylor-Abel-Poisson linear summation method of the Fourier series are   considered  for functions of several variables, periodic with respect to the hexagonal domain. This type of periodicity is  defined by a  lattice which is a discrete subgroup defined by $A{\mathbb Z}^d$, where $A$ is a nonsingular matrix and the periodic function satisfies the relation $f({\bf x}+A{\bf k}) = f({\bf x})$ for all ${\bf k}\in {\mathbb Z}^d$ and ${\bf x}\in {\mathbb R}^d$. With such a periodicity, one works with exponentials or trigonometric functions of the form ${\mathrm e}^{2\pi{\mathrm i}\alpha\cdot {\bf x}}$, where $\alpha$ and $ {\bf x}$ are in certain sets of ${\mathbb R}^d$, not necessarily the usual trigonometric polynomials.
 Lattices appeared prominently in various problems of  analysis and combinatorics (see, for example, \cite{Conway_Sloane_1999}, \cite{Ebeling_1994}). The Fourier series associated with them were studied in the context of information science and
physics, such as sampling theory \cite{Higgins_1996}, \cite{Marks_1996} and multivariate signal processing \cite{Conway_Sloane_1999}.

In the Euclidean plane ${\mathbb R}^2$, besides the standard lattice ${\mathbb Z}^2$ and the rectangular
domain $[-\frac 12,\frac 12)^2$, one of the lattices is the hexagon lattice and the corresponding spectral set is the regular hexagon.
The study of problems connected with Fourier analysis  for functions given on the hexagonal domains goes back to
the papers \cite{Arestov_Berdysheva_2001}, \cite{Sun_2003}.
 Orthogonal polynomials with respect to the area measure on the regular hexagon were studied in \cite{Dunkl_1987}, where an algorithm  was developed for generating an orthogonal polynomial basis.   Let us %also
note the papers \cite{Li_2008},  \cite{Xu_2010} and  \cite{Xu_2016}  in which important
problems of trigonometric approximation and Fourier analysis  on a hexagon  were considered.
In  \cite{Li_2008}, in particular, the discrete Fourier analysis on the regular hexagon is developed in detail.
  \cite{Xu_2010} deals with the problems of the Abel and Ces\`{a}ro summation  of Fourier series, the degree of approximation, and best approximation by trigonometric polynomials  were considered. In \cite{Xu_2010}, especially, the properties of  the Poisson and  Ces\`{a}ro kernels were studied and  direct and inverse theorems were proved in the terms of best approximations of functions by trigonometric polynomials  and their moduli of smoothness.

 Let us also note the papers of Guven (see, for example, \cite{Guven_2013_2}, \cite{Guven_2020} and the references cited therein), in which the author investigated the problems of approximation of functions defined on the hexagonal domains by various linear summation methods of Fourier series.

 In this paper, the approximative properties of the Taylor-Abel-Poisson linear summation method are   considered  in the hexagonal domain. This method defines the operators that possess  the main properties of the Abel--Poisson  and Taylor operators. However, they  can also  be adapted to smoothness properties of functions of arbitrarily large order.
The Taylor-Abel-Poisson operators   $A_{\varrho,r}$  were first studied in \cite{Savchuk_2007} where, in terms of these operators, the author gave the constructive characteristic of Hardy-Lipschitz classes %$H^r_p\mathop{\rm Lip}\alpha$
of functions of one variable, holomorphic on the unit disc of the complex plane. In  \cite{Savchuk_Shidlich_2014},  in
terms of approximation estimates by such operators in the spaces ${\mathcal S}^p$ of Sobolev type, the authors give a constructive description of classes of functions  of several variables whose generalised derivatives belong to the classes ${\mathcal S}^pH_\omega$. Approximations of functions of one variable by similar operators of polynomial type  were studied in \cite{Leis_1963}, \cite{Butzer_Sunouchi_1964}, \cite{Chui_1983}, \cite{Holland_1983}, \cite{Mohapatra_1985}, \cite{Chandra_1988} etc.

 In the integral metrics, for $2\pi$-periodic functions and for functions of  several variables $2\pi$-periodic in each variable, direct and inverse theorems of approximation   by the operators  $A_{\varrho,r}$  were given in the terms of $K$--functionals of functions  generated by their radial derivatives in   \cite{Prestin_Savchuk_Shidlich_2017} and \cite{Prestin_Savchuk_Shidlich_2019} respectively.

In this paper, we prove direct and inverse theorems of approximation by the Taylor-Abel-Poisson operators for functions periodic with respect to hexagonal domains. We also obtain Bernstein-type inequalities for $L_1$-norm of high-order radial derivatives  of  the Poisson kernel.

\section{Preliminaries}

\subsection{Fourier series on the regular hexagon.}\label{Subsetion1.1}

Let us give  basic definitions  of the hexagon lattice and the hexagonal Fourier series. More detailed information
can be found in \cite{Li_2008}, \cite{Xu_2010}, \cite{Sun_2003}.  The hexagonal lattice is given by ${\cal H}{\mathbb Z}^2$, where the matrix ${\cal H}$ and the spectral set $\Omega_H$ are given by
  \[
  {\cal H}=\left(
  \begin{matrix}
  \sqrt{3}\ \ 0\\ -1\ \ 2
  \end{matrix}
  \right),\quad \Omega_{\cal H}=\Big\{(x_1,x_2):\ -1\le x_2,\frac{\sqrt{3}}2x_1\pm\frac 12x_2\le 1\Big\},
  \]
respectively. The reason why $\Omega_{\cal H}$ contains only half of its boundary is described in \cite{Sun_2003}.
We use the homogeneous coordinates $(t_1, t_2, t_3)$ which satisfy the equality $t_1 + t_2 + t_3 = 0$. If we set
\begin{equation}\label{relation_coord}
 t_1:=-
 \frac {x_2}2+\frac{\sqrt{3}x_1}2,\quad t_2:=x_2,\quad t_3:=-
 \frac {x_2}2-\frac{\sqrt{3}x_1}2,
\end{equation}
then  the hexagonal domain $\Omega_{\cal H}$ becomes
 \[
 \Omega=\Big\{(t_1,t_2,t_3) :\quad  t_1 + t_2 + t_3 = 0,\ \ -1\le t_1,t_2,t_3\le  1\Big\}
 \]
which is the intersection of the plane $t_1 + t_2 + t_3 = 0$ with the cube $[-1, 1]^3$ in $\mathbb R^3$.

 Further, for convenience we denote by ${\mathbb R}^3_{\cal H}$ and ${\mathbb Z}^3_{\cal H}$ the sets of all triples
 ${\bf t}=(t_1,t_2,t_3)$ from the corresponding sets  ${\mathbb R}^3$ and ${\mathbb Z}^3$ such that  $t_1+t_2+t_3=0$.

%The relation between $(x_1, x_2)\in \Omega_{\cal H}$ and ${\bf t}\in \Omega$ is given by
%\begin{equation}\label{relation_coord}t_1=-\frac{x_2}2+\frac{\sqrt{3}x_1}2,\quad t_2=x_2,\quad t_3=-\frac{x_2}2-\frac{\sqrt{3}x_1}2. \end{equation}
If we treat $x\in {\mathbb R}^2$ and ${\bf t}\in {\mathbb R}^3_{\cal H}$ as column vectors, then it follows from (\ref{relation_coord}) that
\begin{equation}\label{relation_coord2}
x=\frac 13{\cal H}(t_1-t_3,t_2-t_3)^{\mathrm t \mathrm r}=\frac 13{\cal H}(2t_1+t_2,t_1+2t_2)^{\mathrm t \mathrm r}
\end{equation}
upon using the fact that $t_1 + t_2 + t_3 = 0$. Computing the Jacobian of the change of
variables shows that ${\mathrm d}x=\frac {2\sqrt{3}}3{\mathrm d}t_1{\mathrm d}t_2$.

A function $f$ is called periodic with respect to the hexagonal lattice ${\cal H}$ (or ${\cal H}$-periodic) if
 \[
 f({\bf x})=f({\bf x}+{\cal H}{\bf k}), \quad {\bf k}\in {\mathbb Z}^2,\quad {\bf x}\in {\mathbb R}^2.
 \]
In homogeneous coordinates, ${\bf x}\equiv {\bf y}$ (mod ${\cal H}$) becomes, as easily seen using
(\ref{relation_coord2}), $\mathbf{t}\equiv \mathbf{s}$ (mod 3) defined as
 $%\[
 t_1-s_1\equiv t_2-s_2\equiv t_3-s_3\ (\rm mod\ 3).
 $ %\]
Thus, a function $f(\mathbf{t})$ is ${\cal H}$-periodic, i.e., $f (\mathbf{t}) = f (\mathbf{t}+\mathbf{j})$
whenever $\mathbf{j}\equiv \mathbf{0}$ (mod 3).

Let $L_p=L_p(\Omega)$, $1\le p<\infty$, be the space of all
measurable functions $f$, given on the hexagonal domain $\Omega$, with finite norm
\begin{eqnarray*}
\|f\|_{p}&:=&
\Big(\frac 1{|\Omega_{\cal H}|}\int\limits_{\Omega_{\cal H}} |f(x_1,x_2)|^p{\mathrm d}x_1{\mathrm d}x_2\Big)^{1/p}\\
&=&\Big(\frac 1{|\Omega|}\int\limits_{\Omega} |f({\bf t})|^p
{\mathrm d} {\bf t}\Big)^{1/p} ,
\end{eqnarray*}
where $|\Omega_{\cal H}|$ and $|\Omega|$ denote the areas of $\Omega_{\cal H}$ and $\Omega$ respectively.
As usual, by $L_\infty=L_\infty(\Omega)$ we denote the space of all measurable functions bounded almost everywhere on $\Omega$ with norm
\begin{eqnarray*}
\|f\|_{\infty}&:=&\mathop{\rm ess\,sup}\limits_{(x_1,x_2)\in \Omega_{\cal H}} |f(x_1,x_2)|\\
&=&\mathop{\rm ess\,sup}\limits_{{\bf t}\in \Omega} |f({\bf t})|.
\end{eqnarray*}
The inner product on the hexagonal domain is defined by
\begin{eqnarray}\label{inner_prod}
\langle f,g\rangle _{_{\scriptstyle {\cal H}}}&=&\frac 1{|\Omega_{\cal H}|}\int\limits_{\Omega_{\cal H}} f(x_1,x_2)\overline{g(x_1,x_2)}{\mathrm d}x_1{\mathrm d}x_2\nonumber\\
&=&
 \frac 1{|\Omega|}\int\limits_{\Omega} f({\bf t})\overline{g({\bf t})}{\mathrm d} {\bf t}.
\end{eqnarray}

 For any integer triple ${\bf k}=(k_1,k_2,k_3)\in {\mathbb Z}^3_{\cal H}$ and every point ${\bf t}=(t_1,t_2,t_3)\in {\mathbb R}^3_{\cal H}$ of the plane,  we consider the following trigonometric monomials: %complex function:
\[
 \phi_{\bf k}({\bf t} ):={\rm e}^{\frac{2\pi i}3{\bf k}\cdot {\bf t}},\quad \mbox{\rm where}\ {\bf k}\cdot {\bf t}=k_1t_1+k_2t_2+k_3t_3.
 \]
It is known that these functions are ${\cal H}$-periodic and orthogonal with respect to the  inner product (\ref{inner_prod}) and
for ${\bf k}, {\bf j}\in {\mathbb Z}^3_{\cal H}$, $\langle \phi_{\bf k},\phi_{\bf j}\rangle _{_{\scriptstyle {\cal H}}} =
\delta_{{\bf k},{\bf j}}$. Moreover, the set $\{\phi_{\bf j}:\ {\bf j}\in {\mathbb Z}^3_{\cal H}\}$ is an
orthonormal basis of $L_2(\Omega)$ \cite{Sun_2003},  and for any function $f\in L_1(\Omega)$, the Fourier series with respect to the system $\phi$ has the form
 \begin{equation}\label{Fourier_series_f}
S[f]({\bf t})=%\sim  %S[f]({\bf t})=
\sum_{{\bf k}\in {\mathbb Z}^3_{\cal H}}\widehat {f} ({\bf k}) \phi_{\bf k}({\bf t})
 %=\sum_{{\bf k}\in {\mathbb Z}^3_{\cal H}}\widehat {f} ({\bf k}) {\mathrm e}^{\frac {2\pi i}3(k_1t_2+k_2t_2+k_3t_3)}
 , \quad \ \mbox{\rm where}\ \  \widehat {f} ({\bf k}):=\langle f,\phi_{\bf k}\rangle.
 \end{equation}

%%%%%%%%%%%%%%%%%%%%%%%%%%%%%%%%%%%%%%%%%%%%%%%%%%%%%%%%%%%%%%%%%%%%%%%%%%%%%%%%%%%%%%%%%%%%%%%%%%%%%%%%%%%%%%%%
%\subsection{Poisson integral }\label{Subsetion1.2}

\subsection{Poisson integral  and Taylor-Abel-Poisson means }\label{Subsetion1.3}%\subsection{Poisson integral.}\label{Subsetion1.2}

Let $f\in L_1 (\Omega)$. Set    ${\mathbb J}_\nu:=\{{\bf k}\in {\mathbb Z}^3_{\cal H}: |{\bf k}|:=\max_j\{|k_j|\}=\nu\},\ \nu=0,1,\ldots,$ and for an arbitrary number  $\varrho\in [0,1)$, denote by $P(f)\left(\varrho,{\bf t}\right)$ the Poisson integral of $f$, i.e.,
 \begin{equation}\label{Poisson operator}
P(f)\left(\varrho,{\bf t}\right):=\frac 1{|\Omega|}\int_{\Omega}
f({\bf t}+{\bf s})P(\varrho,{\bf s}){\rm d}{\bf s},
\end{equation}
where $P(\varrho ,{\bf t})$ is the  Poisson kernel, i.e.,
 \begin{equation}\label{Poisson kernel}
P(\varrho ,{\bf s}):=\sum_{\nu=0}^{\infty}\varrho^{\nu}\sum_{{\bf k}\in {\mathbb J}_\nu}\phi_{\bf k}({\bf s}).
\end{equation}
For any  $\varrho\in [0,1)$ and $r\in\mathbb N$, consider the transformation
\begin{equation}\label{def Ar}
A_{\varrho,r}(f)({\bf t}):=%A^\vartriangle_{\varrho,r}(f)({\bf t})=
 \sum_{\nu=0}^{\infty}\lambda_{\nu,r}(\varrho)\sum_{{\bf k}\in {\mathbb J}_\nu}\widehat {f} ({\bf k})\phi_{\bf k}({\bf t}),
\end{equation}
where for $\nu=0,1,\ldots,r-1$, the coefficients are defined by $\lambda_{\nu,r}(\varrho)\equiv 1$ and for $v=r,r+1,\ldots$,
\begin{eqnarray}\label{lambda for H^r}
\lambda_{\nu,r}(\varrho)&:=&\sum_{j=0}^{r-1}
{\nu\choose j}(1-\varrho)^{j}\varrho^{\nu-j}\nonumber\\
&=&
\sum_{j=0}^{r-1}\frac{(1-\varrho)^{j}}{j!}~\frac{\partial^j}{\partial\varrho^j}\varrho^{\nu},\quad
 %r\in\mathbb N,~
\varrho\in[0,1).
\end{eqnarray}
The transformation  $A_{\varrho,r}$ can be considered as a linear operator
on  $L_1(\Omega)$ into itself. Indeed, $\lambda_{k,r}(0){=}0$ and for all  $k=r,r+1,\ldots$ and $\varrho\in(0,1)$, we have
\[
\sum_{j=0}^{r-1}
{\nu\choose j}(1-\varrho)^j\varrho^{\nu-j}\le
rq^{\nu}\nu^{r-1},~\mbox{where}~0<q:=\max\{1-\varrho,\varrho\}<1.
\]
Therefore, for any function $f\in L_1(\Omega)$ and for any $0<\varrho<1$, the series on the right-hand side of (\ref{def Ar}) is majorized by the convergent series $2r\|f\|_1\sum_{\nu=r}^{\infty}q^{\nu}\nu^{r-1}$.

{\lemma\label{Lemma1}   Assume that  $f\in L_1(\Omega)$. Then, for any numbers $r\in \mathbb N,$ $\varrho\in[0,1)$ and ${\bf t}\in  \Omega$,
\begin{equation}\label{A_P}
A_{\varrho, r}(f)({\bf t})=\sum_{k=0}^{r-1} \frac{(1-\varrho)^k}{k!}\cdot
\frac{\partial^k }{\partial\varrho^k}P(f)\left(\varrho,{\bf t}\right).
\end{equation}
}

\noindent {\it Proof.} By virtue of %Combining relations
(\ref{Poisson operator}) and (\ref{Poisson kernel}), we have %obtain
the following decomposition of the Poisson integral into a uniformly convergent series:
\begin{equation}\label{series fo Poisson}
P(f)\left(\varrho,{\bf t}\right)=\sum_{\nu=0}^{\infty}\varrho^{\nu}\sum_{{\bf k}\in {\mathbb J}_\nu}\widehat {f} ({\bf k})\phi_{\bf k}({\bf t}) \quad\mbox{for all}~\varrho\in[0,1),~
{\bf t}\in  \Omega.
\end{equation}
Differentiating (\ref{series fo Poisson})  with respect to  $\varrho$, we realize that for any $k=0,1,\ldots$
\begin{equation}\label{Poisson_Derivative}
\sum_{k=0}^{r-1} \frac{(1-\varrho)^k}{k!}
\frac{\partial^k }{\partial\varrho^k}P(f)\left(\varrho,{\bf t}\right)=\sum_{k=0}^{r-1}\sum_{\nu=k}^{\infty}{\nu\choose k}(1-\varrho)^k
 \varrho^{\nu-k}\sum_{{\bf k}\in {\mathbb J}_\nu}\widehat {f} ({\bf k})\phi_{\bf k}({\bf t}) .
\end{equation}
By changing the summation order on the right-hand side of (\ref{Poisson_Derivative}) and using the identity
\begin{equation}\label{identy}
\sum_{k=0}^{\nu}
{\nu\choose k}(1-\varrho)^k\varrho^{\nu-k}=1,~\nu=0,1,\ldots,
\end{equation}
we obtain
\begin{eqnarray*}
\sum_{k=0}^{r-1}\frac{(1-\varrho)^k}{k!}
\frac{\partial^k }{\partial\varrho^k}P(f)\left(\varrho,{\bf t}\right)&=&
\sum_{\nu=0}^{r-1}\sum_{k=0}^{\nu}{\nu\choose k}(1-\varrho)^k\varrho^{\nu-k} \sum_{{\bf k}\in {\mathbb J}_\nu}\widehat {f} ({\bf k})\phi_{\bf k}({\bf t})\\
&&+\sum_{\nu=r}^{\infty}\sum_{j=0}^{r-1}
{\nu\choose j}(1-\varrho)^{j}\varrho^{\nu-j}\sum_{{\bf k}\in {\mathbb J}_\nu}\widehat {f} ({\bf k})\phi_{\bf k}({\bf t})\\
&=&A_{\varrho, r}(f)({\bf t}).
\end{eqnarray*}
\vskip -3mm$\hfill\Box$

\subsection{Radial derivatives and $K$-functionals }\label{Subsetion2.1}

If for a function $f\in L_1(\Omega)$ and for a positive integer  $n$, there exists a function $g\in L_1(\Omega)$  such that
 \[
\widehat g_{\bf k}=\left\{
     \begin{matrix}
     0,\hfill  & \quad \ \mbox{\rm if}\quad |{\bf k}|<n,\quad \  \\  {\displaystyle \frac{|{\bf k}|!}{(|{\bf k}|-n)!}}\widehat f_{\bf k},\quad   \hfill & \mbox{\rm if}\quad |{\bf k}| \ge n,
     \end{matrix}
     \right.\quad  {\bf k}\in {\mathbb Z}^3_{\cal H},
\]
then we say that for the function $f$, there exists the radial derivative  $g$ of order $n$, for which we use the notation $f^{[n]}$.

%%%%%%%%%%%%%%%%%%%%%%%%%%%%%%%%%%%%%%%%%%%%%%%%%%%%%%%%%%%%%%%%%%%%%%%%%%%%%%%%%%%%%%%%%%%%%%%%%%%%%%%%%%%

Let us note that if the function $f^{[n]}\in L_1(\Omega)$, then its Poisson integral can be represented
%%%in the form
as
\begin{equation}\label{diff f[n]}
P(f^{[n]})(\varrho,{\bf t})= P(f)^{[n]}(\varrho, {\bf t})=\varrho^n\frac{\partial^n
}{\partial\varrho^n}P(f)\left(\varrho,{\bf t}\right),
\quad \varrho\in[0,1),~{\bf t}\in  \Omega.
\end{equation}

%%%%%%%%%%%%%%%%%%%%%%%%%%%%%%%%%%%%%%%%%%%%%%%%%%%%%%%%%%%%%%%%%%%%%%%%%%%%%%%%%%%%%%%%%%%%%%%%%%%%%%%%%%%

In the space $L_p(\Omega)$,
the %%% added
$K$--functional of
%%%the
a
function $f$ (see, for example, \cite[Chap.~6]{DeVore_Lorentz_1993})  generated by the radial derivative of order  $n$,  is the following quantity:
\[
K_n(\delta, f)_p:=\inf\left\{\left\|f-h\right\|_p+\delta^n\left\|h^{[n]}\right\|_p: h^{[n]}\in L_p(\Omega)\right\},\quad\delta>0.
\]

%%%%%%%%%%%%%%%%%%%%%%%%%%%%%%%%%%%%%%%%%%%%%%%%%%%%%%%%%%%%%%%%%%%%%%%%%%%%%%%%%%%%%%%%%%%%%%%%%%%%%%%%%%%%%%%%
\section{Main results}
\subsection{Bernstein-type inequality}\label{Bernstein-type-inequality}

In the following assertion, we  give Bernstein-type inequalities for the $L_1$-norm of radial derivatives  of arbitrary order of the Poisson kernel.

{\theorem\label{Th0} For any $r=0,1,\ldots,$ and for any $0\le\varrho<1$, the following inequality holds:
\begin{equation}\label{Integral_main}
I (\varrho):=\int\limits_{\Omega}\left| \frac{\partial^r }{\partial\varrho^r}P(\varrho, {\bf t})\right| \frac {{\rm d}{\bf t}}{|\Omega|} \le \frac{C_r}{(1- \varrho)^{r}},
\end{equation}
where   the constant $C_r\ge 1$ depends  only on $r$.
}

Before proving Theorem \ref{Th0}, let us give a few auxiliary results.

Let  $P_\varrho(t)$  be the usual Poisson kernel, i.e.,
 \[
 P_\varrho(t)=\frac {1-\varrho^2}{q_\varrho(t)},
 \]
 where $\varrho\in[0,1)$,  $t\in {\mathbb R}$ and $q_\varrho(t) = 1-2\varrho\cos t+\varrho^2$.

{\lemma\label{Lemma003}  For any  $r=0,1,2,\ldots,$ and for any $0\le\varrho<1$, the following relation holds:
\begin{equation}\label{derivative one dimentional_Poisson_Kernal|}
 \left|\frac{\partial^{\,r}}{\partial\varrho^{\,r}} P_\varrho(t)\right|\le \frac{2r!}{(1- \varrho)^{r+1}}.
\end{equation}}

\noindent {\it Proof.} Since
\begin{eqnarray}\label{|Classical_Poisson_Kernal|}
P_\varrho(t)&=&\frac{1-\varrho^2}{|1-\varrho\mathrm{e}^{ \mathrm{i}t}|^2}\nonumber\\
&=&\frac{1}{1-\varrho\mathrm{e}^{ \mathrm{i}t} }+
\frac{1}{1-\varrho \mathrm{e}^{-\mathrm{i}t} }-1
,
\end{eqnarray}
we haven for any $r=0,1,2,\ldots$
\[
\frac{\partial^{\,r}}{\partial\varrho^{\,r}} P_\varrho(t)=
\frac{r!\, \mathrm{e}^{ \mathrm{i}rt}}{(1-\varrho\mathrm{e}^{ \mathrm{i}t})^{r+1}}+
\frac{r!\,\mathrm{e}^{- \mathrm{i}r t}}{(1-\varrho\mathrm{e}^{- \mathrm{i}t})^{r+1}},
 \]%
and the inequality  (\ref{derivative one dimentional_Poisson_Kernal|}) follows.
 \vskip -3mm$\hfill\Box$

{\lemma\label{Lemma03}  Assume that
$z_1=\frac{2\pi}{3}(t_2-t_3)$, $z_2=\frac{2\pi}{3}(t_3-t_1)$, $z_3=\frac{2\pi}{3}(t_1-t_2)$, $r=0,1,\ldots,$ and $0\le\varrho<1$.
 The following assertions holds:

i) for any $j=1,2,3$,
\begin{equation}\label{I_1}
I_{1,r}(\varrho):=\int\limits_{\Omega}
 \left|\frac{\partial^{\,r}}{\partial\varrho^{\,r}}P_\varrho( z_j)\right|\frac {{\rm d}{\bf t}}{|\Omega|} \le \frac{2r!}{(1- \varrho)^{r}}.
\end{equation}

ii) for any numbers $j, k=1,2,3$  and  $r_1,r_2=0,1,\ldots$ such that $j\not =k$ and $r_1+r_2=r$,
\begin{equation}\label{I_2}
I_{2,r}(\varrho):=\int\limits_{\Omega} \left|\frac{\partial^{\,r_1}}{\partial\varrho^{\, r_1}}P_\varrho( z_j)
\frac{\partial^{\,r_2}}{\partial\varrho^{\, r_2}}P_\varrho( z_k)\right|
\frac {{\rm d}{\bf t}}{|\Omega|}\le \frac{4 r_1!r_2!}{(1- \varrho)^{r}}.
\end{equation}

iii) for any  $r_j=0,1,\ldots$, $j=1,2,3$,  such that $r_1+r_2+r_3=r$,
\begin{eqnarray}\label{I_3}
I_{3,r}(\varrho)&:=&\int\limits_{\Omega}\left|\frac{\partial^{\,r_1}}{\partial\varrho^{\, r_1}}P_\varrho( z_1)\frac{\partial^{\,r_2}}{\partial\varrho^{\, r_2}}P_\varrho( z_2)\frac{\partial^{\,r_3}}{\partial\varrho^{\, r_1}}P_\varrho( z_3)\right|\frac {{\rm d}{\bf t}}{|\Omega|}\nonumber\\
&\le&\frac{8 r_1!r_2!r_3!}{(1- \varrho)^{r+1}}.
\end{eqnarray}
}

\noindent {\it Proof.} First, consider the case $r=0$. In this case, the integrals in the inequalities (\ref{I_1})--(\ref{I_3}) can be found exactly.  Let us determine the integral $I_{1,0}(\varrho)$. Consider, for example, the case $j=1$ (the proof in the other cases is similar). Then, $z_j=z_1=\frac{2\pi}{3}(t_2-t_3)$ and by virtue of (\ref{|Classical_Poisson_Kernal|}), we have
\begin{eqnarray}\label{Est_I_1_r_0_hp}
I_{1,0}(\varrho)&=&\int\limits_{\Omega} P_\varrho( z_1) \frac {{\rm d}{\bf t}}{|\Omega|}\nonumber\\
&=&
(1-\varrho^2)\int\limits_{\Omega} \frac 1{|1-\varrho{\mathrm e}^{{\mathrm i}z_1}|^2}\frac {{\rm d}{\bf t}}{|\Omega|}\nonumber\\
&=&(1-\varrho^2)\int\limits_{\Omega} \Big|\sum_{l=0}^\infty \varrho^{l}{\mathrm e}^{{\mathrm i}lz_1}\Big|^2\frac {{\rm d}{\bf t}}{|\Omega|}\nonumber\\
&=&(1-\varrho^2)\int\limits_{\Omega} \Big|\sum_{l=0}^\infty \varrho^{l}{\mathrm e}^{\frac {2\pi{\mathrm i}l(t_2-t_3)}3 }\Big|^2\frac {{\rm d}{\bf t}}{|\Omega|}.
\end{eqnarray}
In (\ref{Est_I_1_r_0_hp}), the last sum   does not contain two terms with identical harmonics. Moreover,   each harmonic in it is equal to $\phi_{{\bf k}_l}({\bf t})$, where ${\bf k}_l=(0,l,-l)\in {\mathbb Z}^3_{\cal H}$. Therefore, using  the orthonormality of $\{\phi_{{\bf k}}\}$, we get
\begin{eqnarray*}
I_{1,0}(\varrho)&=&(1-\varrho^2)\sum_{l=0}^\infty \varrho^{2l}\\
&=&(1-\varrho^2)\frac 1{ 1-\varrho^2 }\\
&=&1.
\end{eqnarray*}

%%%%%%%%%%%%%%%%%%%%%%%%%%%%%%%%%%%%%%%%%%%%%%%%%%%%%%%%%%%%%%%%%%%%%%%%%%%%%%%%%%%%%%%%%%%%%%%%%%%%%
Let us find the integral $I_{2,0}(\varrho)$. Similarly, consider, for example the case  $j=1$ and $k=2$. Then, $z_j=z_1=\frac{2\pi}{3}(t_2-t_3)$, $z_k=z_2=\frac{2\pi}{3}(t_3-t_1)$, and
\begin{eqnarray}\label{Series}
\frac {I_{2,0}(\varrho)}{(1-\varrho^2)^2}&=&\frac {1}{(1-\varrho^2)^2}\int\limits_{\Omega} P_\varrho( z_1)P_\varrho( z_2) \frac {{\rm d}{\bf t}}{|\Omega|}\nonumber\\
&=&
 \int\limits_{\Omega}
\frac 1{|(1-\varrho{\mathrm e}^{{\mathrm i}z_1})(1-\varrho{\mathrm e}^{{\mathrm i}z_2})|^2}\frac {{\rm d}{\bf t}}{|\Omega|}\nonumber\\
&=&\int\limits_{\Omega} \left|\sum_{l,n=0}^\infty \varrho^{l+n}{\mathrm e}^{{\mathrm i}(lz_1+nz_2)}\right|^2\frac {{\rm d}{\bf t}}{|\Omega|}\nonumber\\
&=&\int\limits_{\Omega} \left|\sum_{l,n=0}^\infty \varrho^{l+n}{\mathrm e}^{\frac{2\pi{\mathrm i}}{3}(-nt_1+lt_2+(n-l)t_3)}\right|^2\frac {{\rm d}{\bf t}}{|\Omega|}.
\end{eqnarray}
It is easy to show that  the last sum  in (\ref{Series})  does not contain two terms with identical harmonics. Moreover,   each harmonic in it is equal to $\phi_{{\bf k}_{l,n}}({\bf t})$, where ${\bf k}_{l,n}=(-n,l,n-l)\in {\mathbb Z}^3_{\cal H}$. Therefore, using  the orthonormality of $\{\phi_{{\bf k}}\}$, we get
\begin{eqnarray}\label{Est_I_2_r_0}
I_{2,0}(\varrho)&=&(1-\varrho^2)^2\sum_{l,n=0}^\infty \varrho^{2(l+n)}
\nonumber\\
&=&(1-\varrho^2)^2
\frac 1{(1-\varrho^2)^2}\nonumber\\
&=&1.
\end{eqnarray}

In \cite[Proposition 3.1]{Xu_2010}, it was shown  that for all ${\bf t}\in  \Omega $, the Poisson kernel $P(\varrho, {\bf t})$ satisfies the following relations:
\begin{eqnarray}\label{Poisson_kernel_representation_XU}
P(\varrho, {\bf t})&=&\frac{(1-\varrho)^3(1-\varrho^3)}{q_\varrho\Big(\frac{2\pi}{3}(t_1-t_2)\Big)
q_\varrho\Big(\frac{2\pi}{3}(t_2-t_3)\Big)q_\varrho\Big(\frac{2\pi}{3}(t_3-t_1)\Big)}\nonumber\\
&&+
\frac{\varrho(1-\varrho)^2}{q_\varrho\Big(\frac{2\pi}{3}(t_1-t_2)\Big)
q_\varrho\Big(\frac{2\pi}{3}(t_2-t_3)\Big)}\nonumber\\
&&+\frac{\varrho(1-\varrho)^2}{q_\varrho\Big(\frac{2\pi}{3}(t_1-t_2)\Big)
q_\varrho\Big(\frac{2\pi}{3}(t_3-t_1)\Big)}
\nonumber\\
&&+
\frac{\varrho(1-\varrho)^2}{q_\varrho\Big(\frac{2\pi}{3}(t_2-t_3)\Big)q_\varrho\Big(\frac{2\pi}{3}(t_3-t_1)\Big)},
\end{eqnarray}
and
\begin{equation}\label{Poisson_kernel_2}
 \int\limits_{\Omega} P(\varrho, {\bf t}) \frac {{\rm d}{\bf t}}{|\Omega|} =1.
\end{equation}
Taking into account the notations above, relation (\ref{Poisson_kernel_representation_XU}) can be represented as
\begin{eqnarray}\label{Poisson_Kernal_representation_XU_NEW}
P(\varrho, {\bf t})&=&\frac{1-\varrho^3}{(1+\varrho)^3} P_\varrho( z_1)P_\varrho( z_2)P_\varrho( z_3)\nonumber\\
&+&
\frac{\varrho}{ (1+\varrho)^2}\left( P_\varrho( z_1)P_\varrho( z_2)
+P_\varrho(z_1)P_\varrho(z_3)+P_\varrho(z_2)P_\varrho(z_3)\right).
\end{eqnarray}
Combining (\ref{Poisson_Kernal_representation_XU_NEW}), (\ref{Poisson_kernel_2}) and  (\ref{Est_I_2_r_0}), we get
\begin{eqnarray}\label{Est_I_3_r_0}
 I_{3,0}(\varrho)&=&\int\limits_{\Omega}\! P_\varrho( z_1)P_\varrho( z_2)P_\varrho( z_3)
  \frac {{\rm d}{\bf t}}{|\Omega|}\nonumber\\
&=&
 \left(1-\frac {3 \varrho}{(1+\varrho)^2}\right) \frac {(1+\varrho)^3}{1-\varrho^3}\nonumber\\
&=&\frac {1+\varrho^3}{1-\varrho^3}.
\end{eqnarray}

Now let $r=1,2,\ldots$.
By virtue of (\ref{derivative one dimentional_Poisson_Kernal|})
and the Parseval's identity, for any $j=1,2,3$, we have
 \begin{eqnarray}\label{Est_I_1}
I_{1,r}(\varrho)
&\le& \frac{2r!}{(1-\varrho)^{r-1}}
 \int\limits_{\Omega}
  \frac 1{|1-\varrho{\mathrm e}^{{\mathrm i}z_j}|^2}\frac {{\rm d}{\bf t}}{|\Omega|}\nonumber\\
&=&\frac{2r!}{(1-\varrho)^{r-1}}
\frac 1{ 1-\varrho^2 }\nonumber\\
&\le& \frac{2r!}{(1-\varrho)^{r}}.
\end{eqnarray}

If one of the numbers $r_i$, $i=1,2$, is zero, for example, $r_1=0$, then $r_2=r$ and taking into account (\ref{|Classical_Poisson_Kernal|}) and (\ref{derivative one dimentional_Poisson_Kernal|}),  we have
 %%%%%%%%%%%%%%%%%%%%%%%%%%%%%%%%%%%%%%%%%%%%%%%%%%%%%%%%%%%%%%%%%%%%%%%%
 \begin{eqnarray}\nonumber
  I_{2,r}(\varrho)  &=&
\int\limits_{\Omega}  P_\varrho( z_j)
\left|\frac{\partial^{\,r}}{\partial\varrho^{\, r}}P_\varrho( z_k)\right| \frac {{\rm d}{\bf t}}{|\Omega|}
 \\ \nonumber
&\le&   \frac{2r!(1-\varrho^2)  }{(1-\varrho)^{r-1}}
  \int\limits_{\Omega}\!
  \frac 1{|(1-\varrho{\mathrm e}^{{\mathrm i}z_j})(1-\varrho{\mathrm e}^{{\mathrm i}z_k})|^2}\frac {{\rm d}{\bf t}}{|\Omega|}
 \\ \label{Est_I_2}
&\le &   \frac{2r!(1-\varrho^2)  }{(1-\varrho)^{r-1}}   \frac 1{(1-\varrho^2)^2}\nonumber\\
&\le& \frac{2r!  }{(1-\varrho)^{r}} .
\end{eqnarray}
  %%%%%%%%%%%%%%%%%%%%%%%%%%%%%%%%%%%%%%%%%%%%%%%%%%%%%%%%%%%%%%%%%%%%%%%%
If $r_1$, $r_2\not=0$,  then by virtue of (\ref{derivative one dimentional_Poisson_Kernal|}),
we similarly get the estimate (\ref{I_2}):
\begin{eqnarray*}
I_{2,r}(\varrho)&\le&
   \frac{4 r_1!r_2!}{(1-\varrho)^{r_1+r_2-2}}
  \int\limits_{\Omega}
  \frac 1{|(1-\varrho{\mathrm e}^{{\mathrm i}z_j})(1-\varrho{\mathrm e}^{{\mathrm i}z_k})|^2}\frac {{\rm d}{\bf t}}{|\Omega|}\\
&\le&\frac{4 r_1!r_2!}{(1-\varrho)^{r}} .
\end{eqnarray*}
In the case when one or two of the $r_j$ are equal to zero, the estimates of the integral $I_{3,r}(\varrho)$ are obtained similarly to (\ref{Est_I_2}).
If all $r_j\not =0$, then the inequality (\ref{I_3}) is easily obtained by applying estimates (\ref{derivative one dimentional_Poisson_Kernal|}) and (\ref{I_2}):
 %%%%%%%%%%%%%%%%%%%%%%%%%%%%%%%%%%%%%%%%%%%%%%%%%%%%%%%%%%%%%%%%%%%%%%%%
 \begin{eqnarray*}
 I_{3,r}(\varrho)&\le&
\frac{2r_1!}{(1- \varrho)^{r_1+1}}\int\limits_{\Omega}
  \left|\frac{\partial^{\,r_2}}{\partial\varrho^{\, r_2}}P_\varrho( z_2)
  \frac{\partial^{\,r_3}}{\partial\varrho^{\, r_3}}P_\varrho( z_3)\right|
  \frac {{\rm d}{\bf t}}{|\Omega|}
 \\
&\le& \frac{2r_1!}{(1- \varrho)^{r_1+1}}\frac{4r_2!r_3!}{(1- \varrho)^{r_2+r_3}}\\
&=&\frac{8 r_1!r_2!r_3!}{(1- \varrho)^{r+1}}.
\end{eqnarray*}
  %%%%%%%%%%%%%%%%%%%%%%%%%%%%%%%%%%%%%%%%%%%%%%%%%%%%%%%%%%%%%%%%%%%%%%%%
  \vskip -3mm$\hfill\Box$

\noindent {\it Proof of Theorem \ref{Th0}.}
According to (\ref{Poisson_kernel_2}), we see that the inequality (\ref{Integral_main}) is satisfied  for $r=0$. Let us verify that it is also satisfied  in the case $r=1$. By virtue of (\ref{Poisson_Kernal_representation_XU_NEW}), we have
\begin{eqnarray*}
\frac{\partial }{\partial\varrho}P(\varrho, {\bf t})&=&\frac{-3(1+\varrho^2)}{(1+\varrho)^4} P_\varrho( z_1)P_\varrho( z_2)P_\varrho( z_3)\\
&&+\frac{1-\varrho^3}{(1+\varrho)^3} \left(\frac{\partial P_\varrho( z_1)}{\partial\varrho}   P_\varrho( z_2)P_\varrho( z_3)\right.\\
&&+\left.\frac{\partial P_\varrho( z_2)}{\partial\varrho}   P_\varrho( z_1)P_\varrho( z_3) +   \frac{\partial P_\varrho( z_2)}{\partial\varrho}P_\varrho( z_1)   P_\varrho( z_3)
\right)\\
&&+\frac{1-\varrho}{(1+\varrho)^3} \left(P_\varrho( z_1)P_\varrho( z_2) + P_\varrho( z_2)P_\varrho( z_3)+P_\varrho( z_1) P_\varrho( z_3)\right)\\
&&+\frac{\varrho}{(1+\varrho)^2}\left(  \frac{\partial P_\varrho( z_1)}{\partial\varrho}   P_\varrho( z_2)+  \frac{\partial P_\varrho( z_2)}{\partial\varrho}   P_\varrho( z_1)\right.\\
&& +\frac{\partial P_\varrho( z_2)}{\partial\varrho} P_\varrho( z_3)+\frac{\partial P_\varrho( z_3)}{\partial\varrho}   P_\varrho( z_2)\\
&&+\left. \frac{\partial P_\varrho( z_1)}{\partial\varrho}  P_\varrho( z_3) +
  \frac{\partial P_\varrho( z_3)}{\partial\varrho}P_\varrho( z_1) \right).
\end{eqnarray*}
Then, using the estimates (\ref{Est_I_3_r_0}), (\ref{I_3}),  (\ref{Est_I_2_r_0}), (\ref{Est_I_2}),  we get
\begin{eqnarray*}
 \int\limits_{\Omega}\left|\frac{\partial }{\partial\varrho}P(\varrho, {\bf t})\right|\frac {{\rm d}{\bf t}}{|\Omega|}&\le& \frac{3(1+\varrho^2)}{(1+\varrho)^4} \frac {1+\varrho^3}{1-\varrho^3}+\frac{1-\varrho^3}{(1+\varrho)^3} \frac{24  }{(1- \varrho)^{2}}\\
 &&+\frac{3(1-\varrho)}{(1+\varrho)^3}
 + \frac{\varrho}{(1+\varrho)^2} \frac{12  }{1-\varrho}\\
&\le&\frac {C_1}{1-\varrho} .
\end{eqnarray*}
Therefore, in the case $r=1$, the inequality (\ref{Integral_main}) is indeed satisfied. For $r>1$, the validity of (\ref{Integral_main}) is verified in a similar way using the Leibniz rule of differentiation and Lemma \ref{Lemma03}.

\vskip -3mm$\hfill\Box$

{\corollary\label{Lemma3}  Assume that  $f\in L_p(\Omega)$,  $1\le p\le \infty$  and $\varrho\in (0,1)$. Then, for any $r=0,1,\ldots$,
\begin{equation}\label{|derivative_Poisson_Integral|}
\left\| \frac{\partial^r }{\partial\varrho^r}P(f)(\varrho, \cdot)\right\|_{p}\le C_r\frac{ \|f\|_{_{\scriptstyle p}}}{(1- \varrho)^{r}},
\end{equation}
where the constant $C_r$ depends only on $r$.
}

{\it Proof. } For any $f\in L_p(\Omega)$, by virtue of  the integral Minkowski inequality,  we have
 \[% \begin{equation}\label{Poisson kernel_NEWEST0}
\left\| \frac{\partial^r }{\partial\varrho^r}P(f)\left(\varrho,\cdot\right)\right\|_{p} \le \|f\|_{_{\scriptstyle p}}
\int_\Omega\left|\frac{\partial^r}{\partial\varrho ^r}P(\varrho ,{\bf t})\right| \frac {{\rm d}{\bf t}}{|\Omega|}.
\]%\end{equation}
Therefore, to prove Corollary \ref{Lemma3}, it is sufficient to apply relation (\ref{Integral_main}) to estimate the integral on the right-hand side of the last inequality.

 \vskip -3mm$\hfill\Box$

\subsection{Direct and inverse  approximation  theorems }\label{Subsetion2.2}
Let ${\mathscr Z}$ and ${\mathscr Z}_n$, $n\in {\mathbb N}$, denote the sets of all continuous strictly  increasing functions $\omega(t)$, $t\in [0,1]$, with  $\omega(0)=0$ satisfying the following conditions (\ref{Z}) and (\ref{Z_n}), respectively:
\begin{equation}\label{Z}
 \quad\int_0^\delta\frac{\omega(t)}{t}{\mathrm d}t={\mathcal O}(\omega(\delta)),\quad \delta\to 0+
\end{equation}
and
\begin{equation}\label{Z_n}
 \int_\delta^1\frac{\omega(t)}{t^{n+1}}{\mathrm d}t={\mathcal O}\left(\frac{\omega(\delta)}{\delta^n}\right),\quad\delta\to 0+.
\end{equation}
Conditions (\ref{Z}) and (\ref{Z_n})  are %the
well-known %Zygmund--Bari--Stechkin conditions
(see, for example, \cite{Bari_Stechkin_1956}).

{\theorem\label{Th1} Assume that $f\in L_p(\Omega)$, $1\le p\le\infty$, $n, r\in\mathbb N$, $n\le r$ and $\omega\in  {\mathscr Z}$. If there exists the derivative $f^{[r-n]}\in L_p(\Omega)$ and
\begin{equation}\label{K-funct est}
K_{n}\left( {\delta}, f^{[r-n]}\right)_p={\mathcal O}(\omega(\delta)),\quad\delta\to 0+,
\end{equation}
then
\begin{equation}\label{f-Ap est}
\|f-A_{\varrho, r}(f)\|_p={\mathcal O}\left((1-\varrho)^{r-n}\omega(1-\varrho)\right),\quad\varrho\to 1-.
\end{equation}
}

{\theorem\label{Th2}  Assume that $f\in L_p(\Omega)$, $1\le p\le\infty$, $n, r\in\mathbb N$, $n\le r$ and $\omega\in  {\mathscr Z}\cap {\mathscr Z}_{ n}$.
If the relation $(\ref{f-Ap est})$ holds, then $f^{[r-n]}\in L_p(\Omega)$ and  $(\ref{K-funct est})$ also  holds.
}

The proof of the Theorems \ref{Th1} and \ref{Th2} will be given in the Subsection \ref{Proof_of TH2,3}. Here, we give some comments.

{\remark  For a given $n\in\mathbb N$, from condition $(\ref{Z_n})$  it follows that
 \[
 \mathop{\rm lim~inf}\limits_{\delta\to 0+}(\delta^{-n}\omega(\delta))>0
 \]
 or, equivalently,
 that $(1-\varrho)^{r-n}\omega(1-\varrho)\le C (1-\varrho)^r$ as $\varrho\to~1-$. Therefore, if
 condition $(\ref{Z_n})$  is satisfied, then the quantity on the right-hand side of $(\ref{f-Ap est})$ decreases to zero as $\varrho{\to} 1-$ not faster than the function $(1-\varrho)^r$. }

Here and below, $C$,  $C_{r,n}$, $c_r$, $c_{r,n}, ...$ are positive numbers that do not depend on $\varrho$.

Consider the case %where
$\omega(t)=t^\alpha$, $\alpha>0$.

{\corollary\label{Cor1} Assume that $f\in L_p(\Omega)$, $1\le p\le\infty$, $n, r\in\mathbb N$, $0<\alpha\le n\le r$. The following two relations are equivalent:
\begin{equation}\label{f-Ap est-alpha}
\|f-A_{\varrho, r}(f)\|_p={\mathcal O}\left((1-\varrho)^{r-n+\alpha}\right),\quad\varrho\to 1-,
\end{equation}
and
\begin{equation}\label{K-funct est-alpha}
f^{[r-n]}\in L_p(\Omega)\quad and\quad  K_{n}\left( {\delta}, f^{[r-n]}\right)_p={\mathcal O}(\delta^\alpha),\quad\delta\to 0+.
\end{equation}

}

{\remark\label{Rem2}  If % %n the case where
$\alpha=n$, then the relation (\ref{f-Ap est-alpha}) does not depend on $n$ and $\alpha$. Therefore, the sets of functions $f\in L_p(\Omega)$ such that  $f^{[r-n]}\in L_p(\Omega)$  and $K_{n}\left( {\delta}, f^{[r-n]}\right)_p={\mathcal O}(\delta^n),\ \delta\to 0+$, coincide for any positive integer $n\le r$.
}

Recall (see, e.g., \cite[Ch.~5]{Butzer_Nessel_1971}, \cite[Ch.~2]{Stepanets_2005M}) that a summation method generated by the operator
$A_{\varrho, r}(f)$ is saturated in the space $ L_p(\Omega)$ if there exists a positive function $\varphi$ (which is called the saturation order) defined on the interval
$[0,1)$, monotonically decreasing to zero as $\varrho\to 1-$, and such that each function $f\in  L_p(\Omega)$ satisfying the relation %for which
\begin{equation}\label{Saturation_1}
\|f-A_{\varrho, r}(f)\|_p=\mbox{\tiny  $\mathcal O$} \left(\varphi(\varrho)\right),\quad \varrho\to 1-,
\end{equation}
is an invariant element of the operator $A_{\varrho, r}(f)$  (i.e., $A_{\varrho, r}(f)=f$) and, furthermore, the set
\begin{equation}\label{Saturation_2}
\Phi(A_{\varrho, r})_{p}=\Big\{f\in L_p(\Omega):\quad \|f-A_{\varrho, r}(f)\|_p={\mathcal O} \left(\varphi(\varrho)\right),\quad \varrho\to 1- \Big\}
\end{equation}
(which is   called the saturation class) contains at least one  noninvariant element.

Setting   $n=\alpha=1$ in Corollary \ref{Cor1}, we see that for the function $f\in L_p(\Omega)$, the relation  \[
\|f-A_{\varrho, r}(f)\|_p={\mathcal O} \left((1-\varrho)^r\right),\quad \varrho\to 1-,
\]
holds, if and only if % the following relation holds:
\begin{equation}\label{Saturation_21}
 K_{1}\left( {\delta}, f^{[r-1]}\right)_p={\mathcal O}(\delta),\quad\delta\to 0+.
 \end{equation}
 Furthermore, by virtue of $(\ref{def Ar})$,  for   ${\bf k}\in {\mathbb J}_\nu$, $\nu=0,1,\ldots$, we have
\[
 \frac 1{|\Omega|}\int_{\Omega} (f({\bf t})-A_{\varrho,r}(f)({\bf t}))\overline{\phi_{\bf k}({\bf t})}{\mathrm d} {\bf t}=\left(1-\lambda_{\nu,r}(\varrho)\right) \widehat {f} ({\bf k}),
\]
where  the coefficients  $\lambda_{\nu,r}(\varrho)= 1$ when $\nu=0,1,\ldots,r-1$, and for $\nu\ge r$, they are given by %relation
 (\ref{lambda for H^r}). This yields that  for any   $\varrho\in [0,1)$, $r\in\mathbb N$, ${\bf k}\in {\mathbb J}_\nu$ and $\nu\ge r$,
 \begin{eqnarray*}
 \|f-A_{\varrho, r}(f)\|_p&\ge&
 \left|1-\lambda_{\nu,r}(\varrho)\right| \left|\widehat{f}({\bf k})\right|\\
 &=&\left(1-\sum_{j=0}^{r-1}
{\nu\choose j}(1-\varrho)^{j}\varrho^{\nu-j}\right) \left|\widehat{f}({\bf k})\right|.
 \end{eqnarray*}
By virtue of (\ref{identy}), for ${\bf k}\in {\mathbb J}_\nu$ and $\nu\ge r$, we have
 %%%%%%%%%%%%%%%%%%%%%%%%%%%%%%%%%%%%%%%%%%%%%%%%%%%%%%%%%%%%%%%%%%%%%%%%
 \begin{eqnarray}\nonumber
\left|\widehat {f} ({\bf k})\right|{\nu\choose r}  &=&
\lim_{\varrho\to 1}  \frac {\left|\widehat {f} ({\bf k})\right|}{(1-\varrho)^r}\sum_{j=r}^{\nu}
{\nu\choose j}(1-\varrho)^{j}\varrho^{\nu-j}
 \\ \nonumber
&=&\left| \widehat {f} ({\bf k})\right|\lim_{\varrho\to 1}  \frac {\left|1-\lambda_{\nu,r}(\varrho)\right|}{(1-\varrho)^r}
 \\ \nonumber
&\ge&   \lim_{\varrho\to 1}  \frac {\|f-A_{\varrho, r}(f)\|_p}{(1-\varrho)^r}.
\end{eqnarray}

 %%%%%%%%%%%%%%%%%%%%%%%%%%%%%%%%%%%%%%%%%%%%%%%%%%%%%%%%%%%%%%%%%%%%%%%%

Therefore,  the %following
relation  $%\begin{eqnarray}\label{Saturation_1}
\|f-A_{\varrho, r}(f)\|_p=\mbox{\tiny  $\mathcal O$} \left((1-\varrho)^r)\right),\ \varrho\to 1-,$
 only holds  in the   case when $f({\bf t})=\sum_{\nu=0}^{r-1}\sum_{{\bf k}\in {\mathbb J}_\nu}\widehat {f} ({\bf k})\phi_{\bf k}({\bf t})$ is a polynomial of order not exceeding $r-1$, that is  $A_{\varrho, r}(f)=f$.

 Thus, we conclude that the linear summation method generated by the  operator
$A_{\varrho, r}(f)$ is saturated in the space $ L_p(\Omega)$. The   saturation order is  the function $\varphi(\varrho)=(1-\varrho)^r$ and the saturation class $\Phi(A_{\varrho, r})_{p}$ is the set of all functions $f\in L_p(\Omega)$ satisfying  (\ref{Saturation_21}).

In  \cite{Savchuk_2007}, a similar fact was proved in the Hardy spaces  $H_p$ of functions of one variable, holomorphic on the unit disc of the complex plane.

\subsection{Example }\label{Subsetion3.3_NEW}

Let ${\bf k_0}=(k_{0,1},k_{0,2},k_{0,3})\in {\mathbb Z}^3_{\cal H}$ be any   triple  such that
$|{\bf k_0}|:=\max_j\{|k_{0,j}|\}=r$, $r=1,2,\ldots$. Consider the function
$f_0=\phi_{{\bf k}_0}({\bf t})={\rm e}^{\frac{2\pi i}3{\bf k_0}\cdot {\bf t}}$. Then
 \begin{eqnarray*}
 \|f_0-A_{\varrho, r}(f_0)\|_p &=& \|(1-\lambda_{r,r}(\varrho))\phi_{{\bf k}_0}({\bf t})\|_p\\
&=&
 \Bigg\|\Big(1-\sum_{j=0}^{r-1}{r\choose j}(1-\varrho)^{j}\varrho^{r-j}\Big)\phi_{{\bf k}_0}({\bf t})\Bigg\|_p=  (1-\varrho)^{r} . %\le {n\choose r}(1-\varrho)^{r}
\end{eqnarray*}
Let us also  show that for $0<\delta\le \frac 1r$
\begin{equation}\label{Saturation_21222}
 K_{1}\left( {\delta}, f^{[r-1]}_0\right)_p=r\cdot r! \cdot\delta.
 \end{equation}
For this, we first prove the following assertion.

{\lemma\label{Exact value for K-funct lem}
Let $1\le p\le\infty$, $n\in\mathbb N$, ${\bf k}\in\mathbb J_n$ and   $f=a+\phi_{\bf k}$, $a\in\mathbb C$. Then
\begin{equation}\label{Exact value for K-funct}	
K_1(\delta,f)_p=\left\{
     \begin{matrix}\delta n, \hfill  & \quad \
     \mbox{if } \  0\le\delta\le\frac{1}{n},\\
     1,\hfill  & \quad \ \mbox{if } \ \delta\ge\frac{1}{n}.~~~~~
     \end{matrix}\right.
\end{equation}
}

\noindent {\it Proof.} Let $0\le\delta\le\frac{1}{n}$. Since $f^{[1]}=n\phi_{\bf k}$, we have the upper estimate
$K_1(\delta,f)_p\le\delta n\|\phi_{\bf k}\|_p=\delta n$.

To prove the lower estimate for a fixed arbitrary small  $\varepsilon>0$, consider a function  $h$  such that $h^{[1]}\in L_p(\Omega)$ and
 \begin{equation}\label{K+epsilon}
\|f-h\|_p+\delta\left\|h^{[1]}\right\|_p\le K_1(\delta,f)_p+\varepsilon.
\end{equation}
From the formulas
\[
\left(1-\widehat h({\bf k})\right)\phi_{\bf k}({\bf t})=
\frac{1}{|\Omega|}\int_\Omega\left(f-h\right)({\bf s})\phi_{\bf k}({\bf t}-{\bf s})\mathrm d{\bf s}
\]
and
\[
\widehat h({\bf k})\phi_{\bf k}({\bf t})=\frac{1}{n|\Omega|}\int_\Omega h^{[1]}({\bf s})\phi_{\bf k}({\bf t}-{\bf s})\mathrm d{\bf s},
\]
we get,  respectively,
\begin{equation}\label{Phi}
\left|1-\widehat h({\bf k})\right|\le\|f-h\|_p
\end{equation}
and
\begin{equation}\label{P-Phi}
\left|\widehat h(\bf k)\right|\le\frac{1}{n}\left\|h^{[1]}\right\|_p.
\end{equation}
In view of (\ref{K+epsilon}) and the obvious estimate $K_1(\delta,f)_p\le 1$, the last inequality gives us
 \begin{eqnarray*}
\left|\widehat h({\bf k})\right|&\le &\|f-h\|_p+\frac{1}{n}\left\|h^{[1]}\right\|_p\\
&\le & K_1\left(\frac{1}{n},f\right)_p+\varepsilon\\
&\le & 1+\varepsilon.
\end{eqnarray*}

Therefore, for all $\delta\in[0,\frac{1}{n}]$ we get
 \begin{eqnarray*}
\left|1-\widehat h({\bf k})\right|&\ge & 1-\left|\widehat h(\bf k)\right|
 \\ &\ge & \delta n\left(1-\left|\widehat h(\bf k)\right|+\varepsilon\right)-\varepsilon.
\end{eqnarray*}
Combining this inequality with (\ref{K+epsilon})--(\ref{P-Phi}), we obtain
 \begin{eqnarray*}
\delta n &\le & \delta n+\left(\left|1-\widehat h({\bf k})\right|-\delta n\left(1-\left|\widehat h(\bf k)\right|+\varepsilon\right)+\varepsilon\right)\\
 &=& \left|1-\widehat h({\bf k})\right|+\delta n\left|\widehat h(\bf k)\right|+\varepsilon\left(1-\delta n\right)\\	 
&\le & K_1(\delta,f)_p+2\varepsilon,
\end{eqnarray*}
which proves (\ref{Exact value for K-funct}) for $0\le\delta\le\frac{1}{n}$, since $\varepsilon$ is arbitrary small.
	
For $\delta\ge\frac{1}{n}$ the result (\ref{Exact value for K-funct})   follows on observing that the function $\delta\mapsto K_1(\delta,f)$ increases on $\mathbb R_+$ while staying bounded by 1.
\vskip -3mm$\hfill\Box$

To prove relation (\ref{Saturation_21222}) it is sufficient to apply Lemma \ref{Exact value for K-funct lem} for the function $f_0^{[r-1]}= r! \phi_{{\bf k}_0}({\bf t})$.

\subsection{Discussion}\label{Subsetion3.4_NEW} Note that the results of Subsection \ref{Subsetion2.2} %3.2
can be applied to characterize classes of functions given by  $K$-fucntionals and moduli of smoothness.
In particular, they can be used to study the properties of classical moduli of smoothness, their
connection with $K$-functionals generated by radial and usual derivatives.
 The methods and results of Subsection \ref{Bernstein-type-inequality}, in particular, Theorem  \ref{Th0},
Lemma \ref{Lemma03}, and Corollary \ref{Lemma3}, are crucial for studying properties of the Poisson kernel in the spaces $L_p(\Omega)$.  In addition, interpolation methods for approximation  over hexagonal grids are discussed in \cite{Li_2008},  \cite{Xu_2010}. The results now available can be used to examine the approximation orders whether they are best possible. The same applies to the papers \cite{Xu_2016} and \cite{Guven_2013_2}, in which certain summation methods are used similar to the Taylor-Abel-Poisson means studied here.

\subsection{Auxiliary statements }\label{Subsetion2.3}
In the proof of the Theorems \ref{Th1} and \ref{Th2}, we mainly use the scheme from  \cite{Savchuk_2007}, \cite{Prestin_Savchuk_Shidlich_2017} and \cite{Prestin_Savchuk_Shidlich_2019} modifying it in due consideration of the
peculiarities of the spaces $L_p(\Omega)$. First, let us give some auxiliary results. The proofs are similar to those of Lemmas 3.6--3.8 in \cite{Prestin_Savchuk_Shidlich_2019}.

%\newpage

For any function $f\in L_p(\Omega)$,  $1\le p \le \infty$, $0\le\varrho<1$ and $r=0,1,2,\ldots$, we define
\begin{equation}\label{M_p_diff}
M_p(\varrho, f,r):=\varrho^r\left\| \frac{\partial^r }{\partial\varrho^r}P(f)\left(\varrho,\cdot\right)\right\|_{p}=
\left\|P(f)^{[r]}(\varrho,\cdot)\right\|_{p}.
\end{equation}

{\lemma\label{Lemma31}  Assume that  $f\in L_p(\Omega)$,  $1\le p\le \infty$,   and $\varrho\in (0,1)$. Then, for any nonnegative integers $r$ and $n$, the following inequality holds:
\begin{equation}\label{|derivative A|}
\|A_{\varrho, n}^{[r]}(f)\|_p \le  C_{r,n}\frac{\|f\|_p}{(1- \varrho)^{r}}
\end{equation}
where the constant
$C_{r,n}$ depends only on $r$ and $n$.
}

%%%%%%%%%%%%%%%%%%%%%%%%%%%%%%%%%%%%%%%%%%%%%%%%%%%%%%%%%%%%%%%%%%%%%%%%%%%%%%%%%%%%%%%%%%%%%%%%%%%%%%%%%%%%%%%%%%%%%%%%
\noindent {\it Proof.}   According to (\ref{A_P}) and (\ref{diff f[n]}), for nonnegative integers $r\in {\mathbb N}$ we have
\begin{eqnarray*}
A_{\varrho,n}^{[r]}(f)({\bf t})&=&\left(\sum_{k=0}^{n-1}\frac{ P(f)^{[k]}(\varrho,\cdot) }{\varrho^k k!}(1-\varrho)^k\right)^{[r]}({\bf t})\\
&=&\sum_{k=0}^{n-1}\frac{(P(f)^{[k]}(\varrho,\cdot) )^{[r]}({\bf t})}{\varrho^k k!}(1-\varrho)^k.
\end{eqnarray*}

Since for any nonnegative integers $k$ and $r$
\begin{equation}\label{diff f[n]_f[k]}
(P(f)^{[r]}(\varrho,\cdot) )^{[k]}({\bf t})=(P(f)^{[k]}(\varrho,\cdot) )^{[r]}({\bf t}),
\end{equation}
%%%then
we obtain
$$
A_{\varrho,n}^{[r]}(f)({\bf t})=\sum_{k=0}^{n-1}\frac{(P(f)^{[r]}(\varrho,\cdot) )^{[k]}({\bf t})}{\varrho^k k!}(1-\varrho)^k.
$$
This yields
\begin{eqnarray*}
\|A_{\varrho, n}^{[r]}(f)\|_p&\le&\sum_{k=0}^{n-1}\frac{\left\|(P(f)^{[r]}(\varrho,\cdot) )^{[k]}\right\|_p}{\varrho^k k!}(1-\varrho)^k\\
&\le& \sum_{k=0}^{n-1}\left\| \frac{\partial^k }{\partial \varrho^k} P(f)^{[r]}(\varrho,\cdot) \right\|_p\frac{(1-\varrho)^k}{  k!} ,
\end{eqnarray*}
where by virtue of   Corollary \ref{Lemma3} and (\ref{M_p_diff})
\begin{eqnarray}\label{Ab1}
\left\| \frac{\partial^k  }{\partial \varrho^k}P(f)^{[r]}(\varrho,\cdot)\right\|_p&\le& C_k\frac{\left\|P(f)^{[r]}(\varrho,\cdot)\right\|_{_{\scriptstyle p}}}{(1-\varrho)^{k}}\nonumber\\
&=&\frac{ C_k}{(1-\varrho)^{k}} M_p(\varrho,f,r).
\end{eqnarray}
Therefore,
\begin{equation}\label{Ab2}
\left\|A_{\varrho,n}^{[r]}(f)\right\|_{_{\scriptstyle p}}\le\sum_{k=0}^{n-1}\frac{C_k}{k! } M_p(\varrho,f,r)
\end{equation}
 and applying Corollary \ref{Lemma3}, we get  inequality  (\ref{|derivative A|}). \vskip -3mm$\hfill\Box$

%%%%%%%%%%%%%%%%%%%%%%%%%%%%%%%%%%%%%%%%%%%%%%%%%%%%%%%%%%%%%%%%%%%%%%%%%%%%%%%%%%%%%%%%%%%%%%%%%%%%%%%%%%%%%%%%%%%

{\lemma\label{Main_Lemma} Assume that  $f\in L_p(\Omega)$,  $1\le p\le \infty$. Then, for any numbers $n\in \mathbb N$ and $\varrho\in[1/2,1)$,
%%% change of equation style
\begin{eqnarray}\label{Ineq_Lemma}
c_{1,n}(1-\varrho)^{n}M_p\left(\varrho, f,n\right)&\le& K_{n}\left(1-\varrho, f\right)_p\\
&\le& c_{2,n}\left(\|f-A_{\varrho, n}(f)\|_p+(1-\varrho)^{n}M_p\left(\varrho, f,n\right)\right)\nonumber,
\end{eqnarray}
%%%%%%%%%%%%%%%%%%%%%%%%%%%%%%%%%%%%%%%%%%%%%%%%%%%%%%%%%%%%%%%%%%%%%%%%%%%%%%
%%%%%%%%%%%%%%%%%%%%%%%%%%%%%%%%%%%%%%%%%%%%%%%%%%%%%%%%%%%%%%%%%%%%%%%%%%%%%%
where   the constants
$c_{1,n} $ and $c_{2,n} $ depend  only on $n$.
}

\noindent {\it Proof.}   First, let us note that the statement of Lemma \ref{Main_Lemma} is trivial in the case
when $f$ is a polynomial of order not exceeding $n-1$ %, i.e., when $f({\bf t})=\sum_{\nu=0}^{n-1}\sum_{{\bf k}\in {\mathbb J}_\nu}\widehat {f} ({\bf k})\phi_{\bf k}({\bf t})$,
as well as in the case    $\varrho=0$. Therefore, from now on we exclude these two cases.

Let $g$ be a function such that  $g^{[n]}\in L_p(\Omega)$. Using Corollary \ref{Lemma3} and relation (\ref{M_p_diff}), we get
\begin{eqnarray}\nonumber
\left\|\frac{\partial^n }{\partial\varrho^n}P(f)\left(\varrho,\cdot\right)\right\|_p &=&
\left\|\frac{\partial^n }{\partial\varrho^n}P(f-g)\left(\varrho,\cdot\right)+ \frac{\partial^n }{\partial\varrho^n}P(g)
\left(\varrho,\cdot\right)\right\|_p
 \\ \nonumber
&\le&  C_n\frac{\|f-g\|_p}{(1- \varrho
%^\frac 13
)^{n}}+\left\|\frac{\partial^n }{\partial\varrho^n}P(g)\left(\varrho,\cdot\right)\right\|_p.
\end{eqnarray}
Taking into account  (\ref{diff f[n]}) and Corollary \ref{Lemma3}, we conclude that for any $\varrho\in (0,1)$,
%%%%%%%%%%%%%%%%%%%%%%%%%%%%%%%%%%%%%%%%%%%%%%%%%%%%%%%%%%%%%%%%%%%%%%%%%%%%%%
\begin{eqnarray}\nonumber
  (1-\varrho)^{n} M_p\left(\varrho, f,n\right) &\le& C_n\varrho^{n}\|f-g\|_p+(1-\varrho)^{n}\left\|P(g^{[n]})\left(\varrho,\cdot\right)\right\|_p\\ \nonumber
&\le&
 C_n\|f-g\|_p+C_0(1-\varrho)^{n}\left\|g^{[n]}\right\|_p.
\end{eqnarray}

%%%%%%%%%%%%%%%%%%%%%%%%%%%%%%%%%%%%%%%%%%%%%%%%%%%%%%%%%%%%%%%%%%%%%%%%%%%%%%

Considering the infimum over all functions $g$ such that $g^{[n]}\in L_p(\Omega)$,  we imply
\[
c_{1,n}(1-\varrho)^{n}M_p\left(\varrho, f,n\right)\le K_{n}\left(1-\varrho, f\right)_p .
\]

On the other hand, from the definition of the  $K$--functional, it follows that
\begin{equation}\label{K_n}
K_{n}\left(1-\varrho, f\right)_p\le\left\|f-A_{\varrho, n}(f)\right\|_p+(1-\varrho)^{n}\left\|\left(A_{\varrho, n}(f)\right)^{[n]} \right\|_p.
\end{equation}
Using estimate (\ref{Ab2}) with $r=n$,  we obtain the right-hand inequality in (\ref{Ineq_Lemma}).  \vskip -0.5mm$\hfill\Box$

%%%%%%%%%%%%%%%%%%%%%%%%%%%%%%%%%%%%%%%%%%%%%%%%%%%%%%%%%%%%%%%%%%%%%%%%%%%%%%%%%%%%%%%%%%%%%%%%%%%%%%%%%%%%%%%%%%%%%%%%%%%%%%%%%%%%%%%%%

 {\lemma\label{Lemma4}  Assume that $f\in L_p(\Omega)$, $1\le p\le \infty$, $0\le \varrho<1$ and $r=2,3,\ldots$ such that
\begin{equation}\label{f-A(f)Condition}
\int_\varrho^1 \left\|\frac{\partial^r}{\partial\zeta^r}P(f)(\zeta,\cdot)\right\|_{_{\scriptstyle p}}(1-\zeta)^{r-1}{\mathrm d}\zeta<\infty.
\end{equation}
Then for almost all ${\bf t}\in  \Omega$,
\begin{equation}\label{f-A(f)}
f({\bf t})-A_{\varrho, r}(f)({\bf t})=\frac{1}{(r-1)!}\int_\varrho^1(1-\zeta)^{r-1}\cdot \frac{\partial^r
}{\partial\zeta^r}P(f)(\zeta,{\bf t}){\mathrm d}\zeta.
\end{equation}
}

\noindent {\it Proof.} For  fixed $r=2,3,\ldots$ and $0\le \varrho<1$, the integral on the right-hand side of (\ref{f-A(f)}) defines a certain function $F({\bf t})$. By virtue of (\ref{f-A(f)Condition}) and  the integral Minkowski inequality, we conclude that the function  $F$ belongs to the space $L_p(\Omega)$. Let us find the Fourier coefficients of $F$ and compare them with  the Fourier coefficients of the function $G:=f-A_{\varrho, r}(f)$.

Since for any $\nu=r,r+1\ldots$,
\[
\frac 1{(r-1)!(\nu-r)!}\int_\varrho^{\varrho_1}\zeta^{\nu-r}(1-\zeta)^{r-1}{\mathrm d}\zeta
=\sum_{j=0}^{r-1}
\frac{\varrho_1^{\nu-j}(1-\varrho_1)^{j}-
\varrho^{\nu-j}(1-\varrho)^{j}}{j!(\nu-j)!},
\]
%\rightarrow \frac 1{\nu!}-\sum_{j=0}^{r-1}\frac{\varrho^{\nu-j}(1-\varrho)^{j}}{j!\cdot(\nu-j)!}\]
we have in view of  (\ref{Poisson_Derivative})  for a fixed $\varrho_1\in (\varrho,1)$
\begin{eqnarray}\nonumber
\lefteqn{ \frac 1{(r-1)!}\int_\varrho^{\varrho_1}(1-\zeta)^{r-1}\frac{\partial^r}{\partial\zeta^r}P(f)(\zeta,{\bf t}){\mathrm d}\zeta}\\ \nonumber
&=&
\sum_{\nu=r}^{\infty} \sum_{{\bf k}\in {\mathbb J}_\nu}
\frac{ \nu!\widehat f_{\bf k}\phi_{\bf k}({\bf t})}{(r-1)!(\nu-r)!}
\int_\varrho^{\varrho_1} \zeta^{\nu-r}(1-\zeta)^{r-1}{\mathrm d}\zeta
\\ \label{Ab3}
&=&\sum_{\nu=r}^{\infty}\sum_{{\bf k}\in {\mathbb J}_\nu} \widehat f_{\bf k} \phi_{\bf k}({\bf t})
\sum_{j=0}^{r-1}
{\nu\choose j} \Big(\varrho_1^{\nu-j}(1-\varrho_1)^{j}-
\varrho^{\nu-j}(1-\varrho)^{j}\Big).
\end{eqnarray}
Now, if in the equality (\ref{Ab3}), the value $\varrho_1$ tends to $1-$, then  we observe that the Fourier coefficients $\widehat{F}_{{\bf k}}$  of the  function $F$ are equal to zero when ${\bf k}\in {\mathbb J}_\nu$, $\nu<r$ and for ${\bf k}\in {\mathbb J}_\nu$,  $\nu\ge r$,
\begin{eqnarray}\label{Fourier coeff_F}
\widehat{F}_{{\bf k}}&=&\widehat f_{\bf k}\left(1-\sum_{j=0}^{r-1}
{\nu\choose j} (1-\varrho)^j\varrho^{\nu-j}\right)\nonumber\\
&=&\left(1-\lambda_{\nu,r}(\varrho)\right)\widehat f_{\bf k}.
\end{eqnarray}
Therefore, for all ${\bf k}\in {\mathbb Z}^3_{\cal H}$ we have  $\widehat{F}_{{\bf k}}=(1-\lambda_{\nu,r}(\varrho))\widehat f_{\bf k}=\widehat{G}_{\bf k}$. Hence, for almost all ${\bf t}\in  \Omega$, the representation (\ref{f-A(f)}) holds.  \vskip -3mm$\hfill\Box$

\subsection{Proof of the Theorems \ref{Th1} and \ref{Th2}.}\label{Proof_of TH2,3}
  {\it  Proof of Theorem \ref{Th1}.} Assume that the function $f$ is such that $f^{[r-n]}\in L_p(\Omega)$ and relation (\ref{K-funct est}) is satisfied.  Let us apply  the left-hand side inequality
%%%in
of
Lemma \ref{Main_Lemma} to the function $f^{[r-n]}$. In view of  (\ref{diff f[n]}) and (\ref{M_p_diff}), we obtain
\[
c_{1,n}(1-\varrho)^{n}M_p\left(\varrho, f,r\right)\le K_{n}\left(1-\varrho, f^{[r-n]}\right)_p.
\]
This yields
\begin{equation}\label{est Mp}
M_p\left(\varrho, f,r\right)={\mathcal O}(1) {(1-\varrho)^{-n}}{\omega(1-\varrho)}, \quad \varrho\to 1-.
\end{equation}

Using  relations (\ref{Z}), (\ref{M_p_diff}), (\ref{est Mp})  and the integral Minkowski inequality, we obtain the following estimate
\begin{eqnarray}\label{est Mp estim}
 \int_\varrho^1 \left\|\frac{\partial^r}{\partial\zeta^r}P(f)(\zeta,\cdot)\right\|_{_{\scriptstyle p}}(1-\zeta)^{r-1}{\mathrm d}\zeta&\le& \int_\varrho^1 M_p\left(\zeta, f,r\right)\frac{(1-\zeta)^{r-1}}{\zeta^r}{\mathrm d}\zeta\nonumber\\
 &\le&
 {c_1 }(1-\varrho)^{r-n}\int_\varrho^1\frac{\omega(1-\zeta)}{1-\zeta}{\mathrm d}\zeta\nonumber\\
 &=&
 {\mathcal O}\left((1-\varrho)^{r-n}\omega(1-\varrho)\right)
\end{eqnarray}
for all $0\le\varrho<1$.

Therefore, for almost all ${\bf t}\in  \Omega$, relation (\ref{f-A(f)}) holds. Hence, by virtue of  (\ref{f-A(f)}),
using the integral Minkowski inequality and (\ref{est Mp estim}),  we
finally %%% added
get (\ref{f-Ap est}):
 \begin{eqnarray*}
\|f-A_{\varrho, r}(f)\|_{_{\scriptstyle p}}&\le&\frac{1}{(r-1)!}\int_\varrho^1 M_p\left(\zeta, f,r\right)\frac{(1-\zeta)^{r-1}}{\zeta^r}{\mathrm d}\zeta \\ &=&{\mathcal O}\left((1-\varrho)^{r-n}\omega(1-\varrho)\right), \quad \varrho\to 1-.
\end{eqnarray*}
\vskip -3mm$\hfill\Box$

\bigskip
\noindent {\it  Proof of Theorem \ref{Th2}.} First, let us note that for any function $f\in L_p(\Omega)$ and all fixed numbers $r,s\in {\mathbb N}$ and $\varrho\in (0,1)$, we have

\begin{eqnarray}\nonumber
\|A_{\varrho, r}^{[s]}(f)\|_{_{\scriptstyle p}}&=&\left\|
 \sum_{\nu=s}^{\infty}\frac{\nu!}{(\nu-s)!}\lambda_{\nu,r}(\varrho)
\sum_{{\bf k}\in {\mathbb J}_\nu} \widehat f_{\bf k} \phi_{\bf k}({\bf t}) \right\|_{_{\scriptstyle p}}
 \\ \nonumber
&\ \le&   2r\|f\|_{_{\scriptstyle p}}\left(\sum_{\nu=s}^{\max\{r,s\}-1}\frac{\nu!}{(\nu-s)!}+\sum_{\nu\ge \max\{r,s\}}q^{\nu}\nu^{s+r-1}\right)\nonumber\\
&<&\infty,
\end{eqnarray}
where $0<q=\max\{1-\varrho,\varrho\}<1$. In the case where $s\ge r$,  the  sum $\sum_{\nu=s}^{s-1}$ is set equal to zero.

Put $\varrho_k:=1-2^{-k},~k\in{\mathbb N},$ and $A_k:=A_k(f):=A_{\varrho_k,r}(f)$. For any
${\bf t}\in  \Omega$
and $s\in {\mathbb N}$, consider the series
\begin{equation}\label{series1}
A_0^{[s]}(f)({\bf t})+\sum\limits_{k=1}^\infty (A_k^{[s]}(f)({\bf t})-A_{k-1}^{[s]}(f)({\bf t})).
\end{equation}
According to the definition of the operator  $A_{\varrho,r}$, we conclude that for any $\varrho_1, \varrho_2\in[0,1)$ and $r\in {\mathbb N}$,
\[
A_{\varrho_1,r}\left(A_{\varrho_2,r}(f)\right)=A_{\varrho_2,r}\left(A_{\varrho_1,r}(f)\right).
\]
By virtue of Lemma \ref{Lemma31} and relation  (\ref{f-Ap est}), for any  $k\in {\mathbb N}$ and  $s\in {\mathbb N}$, we have
\begin{eqnarray}\nonumber
\left\|A^{[s]}_{k}-A^{[s]}_{k-1}\right\|_{_{\scriptstyle p}}
&=&\left\|A^{[s]}_{k}(f-A_{k-1}(f))-A^{[s]}_{k-1}(f-A_{k}(f))
\right\|_{_{\scriptstyle p}}\\ \nonumber
&\le&\left\|A^{[s]}_{k}(f-A_{k-1}(f))\right\|_{_{\scriptstyle p}}+\left\|A^{[s]}_{k-1}(f-A_{k}(f)) \right\|_{_{\scriptstyle p}}\\ \nonumber
&\le& C_{s,r}\frac{\left\|f-A_{k-1}(f)\right\|_{_{\scriptstyle p}}}{(1-\varrho_k)^{s}}+C_{s,r}\frac{\left\|f-A_{k}(f)\right\|_{_{\scriptstyle p}}}
{(1-\varrho_{k-1})^{s}} \\ \label{OA}
&=&
{\mathcal O}\left(\frac{\omega(1-\varrho_{k-1})}{(1-\varrho_{k})^{s-r+n}}\right)+
{\mathcal O}\left(\frac{\omega(1-\varrho_{k})}{(1-\varrho_{k-1})^{s-r+n}}\right).
\end{eqnarray}
Therefore, if $s\le r-n$, then
\begin{equation}\label{series12}
\left\|A^{[s]}_{k}-A^{[s]}_{k-1}\right\|_{_{\scriptstyle p}}={\mathcal O}\left(\omega(1-\varrho_{k-1})\right)=
{\mathcal O}\left(\omega(2^{1-k})\right),\quad k\to \infty.
\end{equation}
Consider the sum  $\sum_{k=1}^N \omega(2^{1-k})$, $N\in {\mathbb N}$. Taking into account the monotonicity of the function $\omega$ and (\ref{Z}), we observe that for all $N\in {\mathbb N}$,
\begin{eqnarray*}
 \sum_{k=1}^N\omega(2^{1-k})&\le& \omega(1)+\int_1^{N} \omega(2^{1-t}){\mathrm d}t\\
 &=&
 \omega(1)+ \int_{2^{1-N}}^1 \frac{\omega(\tau){\mathrm d}\tau}{\tau\ln 2} \\
 &<&\infty.
\end{eqnarray*}

Combining the last relation and (\ref{series12}), we  conclude that  the series in (\ref{series1}) converges in the
%%%metrics
norm
of the space $L_p(\Omega)$, $1\le p\le\infty$. Hence, by virtue of the Banach--Alaoglu theorem,  there exists the following subsequence if
 $0\le s\le r-n $:
\begin{equation}\label{series123}
S^{[s]}_{N_j}({\bf t})=A_0^{[s]}(f)({\bf t})+\sum\limits_{k=1}^{N_j} (A_k^{[s]}(f)({\bf t})-A_{k-1}^{[s]}(f)({\bf t})),\quad j=1,2,\ldots
\end{equation}
of partial sums of this series converging to a certain function  $g\in L_p(\Omega)$ almost everywhere on
%%%$[0, 2\pi]$
$\Omega$
 as $j\to\infty$.

Let us show that $g=f^{[s]}$. For this, let us find the Fourier coefficients of the function $g$. For any fixed  ${\bf k}\in{\mathbb Z}^3_{\cal H}$ and all $j=1,2,\ldots,$ we have
\[
\widehat g_{\bf k}:=\frac 1{|\Omega|}\int_{\Omega}S^{[s]}_{N_j}({\bf t})\overline{{\phi}_{\bf k}({\bf t})}{\mathrm d} {\bf t}+\frac 1{|\Omega|}
\int_{\Omega}(g({\bf t})-S^{[s]}_{N_j}({\bf t}))\overline{{\phi}_{\bf k}({\bf t})}{\mathrm d} {\bf t}.
\]
Since the sequence $\{S^{[s]}_{N_j}\}_{j=1}^\infty$ converges almost everywhere on
%%%$[0, 2\pi]$
$\Omega$
to the function $g$, %then
the second integral on the right-hand side of the last equality tends to zero as $j\to\infty$. By virtue of (\ref{series123}) and the definition of the radial derivative,
%%%the first integral is equal to zero, when $|k|<s$,
for $|{\bf k}|=\nu<s$ the first integral is equal to zero,
and for all $|{\bf k}|=\nu\ge s$,
 $$
\frac 1{|\Omega|}\int_{\Omega}S^{[s]}_{N_j}({\bf t})\overline{{\phi}_{\bf k}({\bf t})}{\mathrm d} {\bf t}=
 \lambda_{\nu,r}(1-2^{-N_j})\frac{\nu!}{(\nu-s)!}\widehat f_{\bf k} \mathop{\longrightarrow}\limits_{j\to\infty} \frac{\nu!}{(\nu-s)!}\widehat f_{\bf k}.
 $$
Therefore, the equality  $g=f^{[s]}$ holds true. Hence, for the function  $f$ and all $0\le s\le r-n $, there exists  the derivative $f^{[s]}$, and $f^{[s]}\in L_p(\Omega)$.

%%%Further
Now,
let us prove the estimate  (\ref{est Mp}). By virtue of (\ref{M_p_diff}), (\ref{OA}), for any $k\in\mathbb N$ and $\varrho\in (0,1)$, we have
\begin{eqnarray}\nonumber
M_p\left(\varrho, A_{k}-A_{k-1},r\right)
&\le&
 \left\|A^{[r]}_{k}-A^{[r]}_{k-1}\right\|_{_{\scriptstyle p}}\nonumber
\\ \nonumber
&=&{\mathcal O}\left(\frac{\omega(1-\varrho_{k-1})}{(1-\varrho_{k})^n}\right)+{\mathcal O}\left(\frac{\omega(1-\varrho_{k})}{(1-\varrho_{k-1})^n}\right)
\\ \nonumber
&=&{\mathcal O}\left(2^{kn}\omega(2^{-k+1})+2^{(k-1)n}\omega(2^{-k})\right)
\\ \label{deviation for AA}
&=&
{\mathcal O}\left(2^{(k-1)n}\omega\left(2^{-(k-1)}\right)\right),\ \  k\to \infty.
\end{eqnarray}
By virtue of  (\ref{M_p_diff}), (\ref{|derivative_Poisson_Integral|})  and (\ref{f-Ap est}),  for any $r\in {\mathbb N}$ and $\varrho\in (0,1)$, we obtain
\begin{eqnarray*}
M_p\left(\varrho, f-A_{\varrho,r}(f),r\right)&=&{\mathcal O}(1) \frac{\left\|f-A_{\varrho,r}(f)\right\|_{p}}{(1-\varrho)^{r}}\\
&=&{\mathcal O}\left(\frac{\omega(1-\varrho)}{(1-\varrho)^{n}}\right), \quad \varrho\to 1-.
\end{eqnarray*}
Therefore, for  $N\to \infty$,
\begin{eqnarray}\label{deviation for f}
M_p\left(\varrho_{_{\scriptstyle N}}, f-A_{N}(f),r\right)&=&{\mathcal O}\left(\frac{\omega(1-\varrho_{_{\scriptstyle N}})}{(1-\varrho_{_{\scriptstyle N}})^n}\right)\\
&=&{\mathcal O}\left(2^{Nn}\omega(2^{-N})\right).
\end{eqnarray}

Consider the sum $\sum_{k=1}^N2^{(k-1)n}\omega(2^{-(k-1)})$, $N\in {\mathbb N}$.  For any  $\omega\in {\mathscr Z}_n$,
the function $\omega(t)/t^n$ almost decreases on $(0,1]$, i.e., there exists a number $c_2>0$ such that  $\omega(t_1)/t_1^n\ge c_2 \omega(t_2)/t_2^n$ for any $0<t_1<t_2\le 1$  (see, for example \cite{Bari_Stechkin_1956}). Therefore,
\begin{eqnarray}\nonumber
\lefteqn{ \sum_{k=1}^N2^{(k-1)n}\omega(2^{-(k-1)})}\\ \nonumber
&\le& c_2\left( 2^{(N-1)n}\omega(2^{-(N-1)})
 +\int_1^{N} 2^{(t-1)n}\omega(2^{-(t-1)}){\mathrm d}t\right)\\ \nonumber
&\le& c_2\left(2^{(N-1)n}\omega(2^{-(N-1)})+\int_{2^{-N+1}}^1\frac{ \omega(\tau)\,{\mathrm d}\tau}{
  \tau^{n+1}\ln 2}\right)\\ \label{summe}
 &=&{\mathcal O}\Big( 2^{(N-1)n}\omega(2^{-(N-1)})\Big)\nonumber\\
 &=&{\mathcal O}\Big(2^{Nn}\omega(2^{-N})\Big),\quad N\to \infty.
\end{eqnarray}
Putting $\varrho=\varrho_{_{\scriptstyle N}}$ and taking into account the relations (\ref{deviation for AA}), (\ref{deviation for f}), (\ref{summe}) and
$$
A_{0}({\bf t})=S_{r-1}(f)({\bf t})=\sum_{|{\bf k}|\le r-1}\widehat f_{\bf k}\phi_{\bf k}({\bf t}),
$$
we obtain
\begin{eqnarray}\nonumber
M_p\left(\varrho_N,f,r\right)&=&M_p\left(\varrho_N,f-S_{r-1}(f),r\right)
\\ \nonumber
&=& M_p\left(\varrho_N,f-A_{\varrho_{_{\scriptstyle N}}}+\sum\limits_{k=1}^N (A_{k}-A_{k-1}),r\right)\nonumber\\
&=&{\mathcal O}\left(\sum_{k=1}^{N}2^{(k-1)n}\omega(2^{-(k-1)})\right)
\\ \label{M_p(varrho_N)}
&=&  {\mathcal O}\left(2^{Nn}\omega(2^{-N})\right)\nonumber\\
&=&{\mathcal O}\left((1-\varrho_{_{\scriptstyle N}})^{-n}\omega(1-\varrho_{_{\scriptstyle N}})\right),\quad  N\to \infty.
\end{eqnarray}

If the  function  $\omega$ satisfies the condition $({\mathscr Z}_n)$, then
$\sup\limits_{t\in [0,1]}({\omega(2t)}/{\omega(t)})<\infty$  (see, for example \cite{Bari_Stechkin_1956}).  Furthermore, for all $\varrho\in[\varrho_{_{\scriptstyle N-1}},\varrho_{_{\scriptstyle
N}}]$, we have $1-\varrho_{_{\scriptstyle N}}\le
1-\varrho\le2(1-\varrho_{_{\scriptstyle N}})$. Hence,  relation (\ref{M_p(varrho_N)}) yields the estimate (\ref{est Mp}).

Now, applying the right-hand side inequality in Lemma \ref{Main_Lemma} to the function   $f^{[r-n]}$, we get
\begin{equation}\label{est Kn}
K_{n}\left(1-\varrho, f^{[r-n]}\right)_p \le  c_{2,n}\Big( \|f^{[r-n]}-A_{\varrho, n}(f^{[r-n]})\|_{_{\scriptstyle p}}+
(1-\varrho)^nM_p({\varrho},f,r)\Big).
\end{equation}
%%%%%%%%%%%%%%%%%%%%%%%%%%%%%%%%%%%%%%%%%%%%%%%%%%%%%%%%%%%%%
By virtue of (\ref{M_p_diff}) and (\ref{est Mp}), we conclude  that for  $\varrho\in [1/2,1)$,
\begin{eqnarray}\nonumber
\int_\varrho^1 \left\|\frac{\partial^n }{\partial\zeta^n}P(f^{[r-n]})(\zeta,\cdot)\right\|_{_{\scriptstyle p}}(1-\zeta)^{n-1}{\mathrm d}\zeta&=&\int_\varrho^1 \left\|P(f)^{[r]}(\zeta,\cdot)  \right\|_{_{\scriptstyle p}}\frac{(1-\zeta)^{n-1}}{\zeta^n} {\mathrm d}\zeta\\ \nonumber
&=&\int_\varrho^1 M_p\left(\zeta, f,r\right)\frac{(1-\zeta)^{n-1}}{\zeta^n} {\mathrm d}\zeta
\\ \label{111} &\le& c_1
\int_\varrho^1\frac{\omega(1-\zeta)}{1-\zeta}{\mathrm d}\zeta\nonumber\\
&=&{\mathcal O}\left(\omega(1-\varrho)\right).
\end{eqnarray}
Therefore, we can apply  Lemma \ref{Lemma4} to the function   $f^{[r-n]}$. Taking into account (\ref{M_p_diff}), we obtain
\[
f^{[r-n]}({\bf t})-A_{\varrho, n}(f^{[r-n]})({\bf t})=\frac{1}{(n-1)!}\int_\varrho^1 P(f)^{[r]}(\zeta,{\bf t})\frac{(1-\zeta)^{n-1}}{\zeta^n} {\mathrm d}\zeta.
\]
Using   the integral Minkowski inequality and (\ref{111}), we %%%obtain
conclude
\begin{eqnarray}\nonumber
\|f^{[r-n]}-A_{\varrho, n}(f^{[r-n]})\|_{_{\scriptstyle p}}&\le& \frac{1}{(n-1)!}\int_\varrho^1 M_p\left(\zeta, f,r\right)\frac{(1-\zeta)^{n-1}}{\zeta^n} {\mathrm d}\zeta\\ \label{112}
&=& {\mathcal O}\left(\omega(1-\varrho)\right), \quad \varrho\to 1-.
\end{eqnarray}
Combining the relations (\ref{est Kn}), (\ref{est Mp}) and (\ref{112}),  we finally get (\ref{K-funct est}).

%%%%%%%%%%%%%%%%%%%%%%%%%%%%%%%%%%%%%%%%%%%%%%%%%%%%%%%%%%%%%%%%%%%%%%%%%%%%%%%%%%%%%%%%%%%%%%%%%%%%%%%%%%%%%%%%%%%%%%%%%%%%%%%%%%%%%%%%%%%%%%%%%%%%%%%%%
\section{Acknowledgments}
This research is partially  supported by the Volkswagen Foundation project ``From Modeling and Analysis to Approximation'' and  the Grant H2020-MSCA-RISE-2019, project number  873071 (SOMPATY: Spectral Optimization: From Mathematics to Physics and Advanced Technology).





\bibliographystyle{elsarticle-num.bst}
%\bibliography{Prestin_Savchuk_Shidlich}

%\end{document}
%\endinput
%%
%% End of file `elsarticle-template-num.tex'.

\end{document}